\documentclass[onecolumn]{article}
\usepackage{amsmath,amssymb,graphicx}
\usepackage{algorithm,algorithmic,setspace}
\usepackage{cite,url, caption,subcaption}

\def\R{\mathbb{R}}

\title{A Short Note on Improved ROSETA}
\author{Hassan Mansour\\
Mitsubishi Electric Research Laboratories\\
	Cambridge, MA 02139\\
	mansour@merl.com}

\begin{document}
%
\maketitle
%

%
\section{Problem Formulation}\label{sec:ProblemFormulation}
We consider the problem of identifying at every time $t$ an $r$-dimensional subspace $\mathcal{U}_t$ in $\R^n$ with $r \ll n$ that is spanned by the columns of a rank-$r$ matrix $U_t \in \R^{n\times r}$ from incomplete and noisy measurements
\begin{equation}\label{eq:bt}
	b_t = \Omega_t(U_t a_t + s_t), 
\end{equation}
where $\Omega_t$ is a selection operator that specifies the observable subset of entries at time $t$, $a_t \in \R^{r}$ are the coefficients specifying the linear combination of the columns of $U_t$, and $s_t \in \R^n$ is a sparse outlier vector.

When the subspace $\mathcal{U}_t$ is stationary, we drop the subscript $t$ from $U_t$ and the problem reduces to robust matrix completion or robust principal component analysis where the task is to separate a matrix $B \in \R^{n \times m}$ into a low rank component $UA$ and a sparse component $S$ using incomplete observations 
$$
	B_{\Omega} = \Omega(UA + S).
$$
Here the columns of the matrices $A$ and $S$ are respectively the vectors $a_t$ and $s_t$ stacked horizontally for all $t \in \{1\dots m\}$, and the operator $\Omega$ specifies the observable entries for the entire matrix $B$.

\section{Improved Robust Online Subspace Estimation and Tracking}\label{sec:ROSETA}
We describe in this section an improved version of the robust online subspace estimation and tracking algorithm (ROSETA) first published in~\cite{ROSETA:2015}. 

We consider the combined loss function $\mathcal{L}(U_{t}, s_t, a_t, e_t)$ that includes an $\ell_2$ data misfit term and a one-norm regularizer applied to the outliers $s_t$
\begin{equation}\label{eq:ROSETAloss}
	\mathcal{L}(U_{t}, s_t, a_t, e_t) = \frac{\mu}{2}\|b_t - (U_{t} a_t + s_t + e_t)\|_2^2 + \lambda\|s_t\|_1 
\end{equation}
where $e_t$ is supported on the complement of $\Omega_t$, hereby denoted $\Omega_t^c$, such that $\Omega_t(e_t) = 0$ and $\Omega_t^c(e_t) = -\Omega_t^c(U_t a_t)$.

\subsection{Update the subspace coefficients}
Following the ROSETA framework, we first fix the subspace matrix $U$ to be the previous estimate $U_{t-1}$ and update the triplet $(s_t, a_t, e_t)$ by solving the LASSO problem
\begin{equation}\label{eq:robustSubspaceFit}
	(s_t, a_t, e_t) = \arg\min\limits_{s, a, e, y} \mathcal{L}(U_{t-1}, s, a, e).
\end{equation}
The solution to~\eqref{eq:robustSubspaceFit} is obtained by applying the following sequence of updates:
\begin{equation}\label{eq:LASSOsolution}
\begin{array}{lll}
	a_t^k &=& U_{t-1}^{\dagger}\left(b_t - s_{t}^{k-1} - e_{t}^{k-1} \right)\\
	e_t^k &=& -\Omega_t^c\left(U_{t-1}a_t^k\right)\\
	s_t^k &=& \mathcal{S}_{\lambda}\left(b_t - U_{t-1}a_t^k - e_t^{k} \right)\\
\end{array}
\end{equation}
where $\mathcal{S}_{\tau}(x) = \mathrm{sign}(x)\cdot\max\{|x| - \tau,0\}$ denotes the element-wise soft thresholding operator with threshold $\tau$, $k$ indicates the iteration number, and $^\dagger$ is the Moore-Penrose pseudo-inverse of a matrix.

\subsection{Update the subspace matrix}
The subspace matrix $U_t$ is then updated by minimizing with respect to $U$ the quadratic approximation $Q(U)$ of $\mathcal{L}(U, s_t, a_t, e_t)$ around the estimate $U_{t-1}$, i.e. we define the quadratic function
\begin{equation}\label{eq:quadratic}
	Q(U) := \mathcal{L}(U_{t-1}) + \nabla \mathcal{L}(U_{t-1})^T (U - U_{t-1}) + \frac{\mu}{2}\mathrm{Tr}\left[(U - U_{t-1})(I_r + a_t a_t^T)(U - U_{t-1})^T\right],
\end{equation}
where $\mu$ is the step-size parameter, and $\mathrm{Tr}[\cdot]$ is the trace function. Consequently, the update equation of the subspace matrix $U$ is given by
\begin{equation}\label{eq:Ut_update}
\begin{array}{ll}
	U_t &= \arg\min\limits_{U} Q(U) \\
		&= U_{t-1} + \frac{1}{\mu}(b_t - (U_{t-1} a_t + s_t  + e_t))a_t^T(I_r + a_t a_t^T)^{-1}
\end{array}
\end{equation}

\subsection{Adaptive parameter selection}
Inspired by the adaptive step size selection in GRASTA~\cite{GRASTA:2012}, we developed a corresponding adaptive step-size parameter for ROSETA. 

The parameter $\mu_t$ controls the speed of convergence of the subspace estimate. In particular, a smaller value of $\mu$ allows for faster adaptability of $U_t$ to a changing subspace (larger descent step), whereas a larger value of $\mu$ only permits a small variation in $U_t$. Consider the descent direction 
\begin{equation}\label{eq:descent}
	D_t = (b_t - (U a_t + s_t  + e_t))a_t^T(1 + a_t^T a_t)^{-1}.
\end{equation}
The parameter $\mu_t$ can then be updated according to
\begin{equation}\label{eq:mu_t}
	\mu_t = \frac{C}{1 + \eta_t},
\end{equation}
where $\eta_t = \min\{\eta_{\max}, \max\{C, \eta_{t-1} + \mathrm{sigmoid}\left(\frac{\langle D_{t-1}, D_{t}\rangle}{\|D_{t-1}\|_F\|D_t\|_F}\right)\}\}$, where $\eta_{\max}$ is a control parameter. Here $\mathrm{sigmoid}(x) = f + 2f/(1+e^{10x})$, for some predefined $f$. 

Similar to GRASTA, the intuition behind choosing such an update rule comes from the idea that if two consecutive subspace updates $D_{t-1}$ and $D_t$ have the same direction, i.e. $\langle D_{t-1}, D_{t}\rangle > 0$, then the target subspace is still far from the current subspace estimate. Consequently, the new $\mu_t$ should be smaller to allow for fast adaptability which is achieved by increasing $\eta_t$.  Similarly, when $\langle D_{t-1}, D_{t}\rangle < 0$, the subspace update seems to bounce around the target subspace and hence a larger $\mu_t$ is needed. The new ROSETA algorithm is summarized in Algorithm \ref{alg:ROSETA}.
\begin{algorithm}[ht]\caption{Robust Subspace Estimation and Tracking}
\begin{spacing}{1}
\begin{algorithmic}[1]\label{alg:ROSETA}
\STATE \textbf{Input} Sequence of measurements $\{b_t\}$, $\eta_{\mathrm{LOW}}$, $\eta_{\mathrm{HIGH}}$ 
\STATE \textbf{Output} Sequences $\{U_t\}$, $\{a_t\}$, $\{s_t\}$
\STATE \textbf{Initialize} $U_0$, $\mu_0$, $\eta_0$, $\eta_{\max}$
\FOR{$t = 1 \dots N$}
\STATE \textbf{Solve for $a_t$, $s_t$, and $e_t$:}
\WHILE{not converged} 
\STATE $a_t^k = U_{t-1}^{\dagger}\left(b_t - s_{t}^{k-1} - e_{t}^{k-1} \right)$
\STATE $e_t^k = -\Omega_t^c\left(U_{t-1}a_t^k\right)$
\STATE $s_t^k = \mathcal{S}_{\lambda}\left(b_t - U_{t-1}a_t^k - e_t^{k} \right)$
\ENDWHILE
\STATE \textbf{Update subspace estimate:}
\STATE \ \ $D_t = (b_t - (U a_t + s_t  + e_t))a_t^T(1 + a_t^T a_t)^{-1}$
\STATE \ \ $U_t = U_{t-1} + \frac{1}{\mu}D_t$
\STATE \textbf{Update parameter $\mu_t$:}
\STATE \ \ $\eta_t = \min\{\eta_{\max}, \max\{C, \eta_{t-1} + \mathrm{sigmoid}\left(\frac{\langle D_{t-1}, D_{t}\rangle}{\|D_{t-1}\|_F\|D_t\|_F}\right)\}\}$
\STATE \ \ $\mu_t = \frac{C}{1 + \eta_t}$
\ENDFOR
\end{algorithmic}
\end{spacing}
\end{algorithm}\vspace{-0.2in}

\bibliographystyle{IEEEbib}
\bibliography{ROSETA}

\begin{thebibliography}{1}

\bibitem{ROSETA:2015}
H.~Mansour and X.~Jiang,
\newblock ``A robust online subspace estimation and tracking algorithm,''
\newblock in {\em 2015 IEEE International Conference on Acoustics, Speech and
  Signal Processing (ICASSP)}, April 2015, pp. 4065--4069.

\bibitem{GRASTA:2012}
J.~He, L.~Balzano, and J.~C.~S. Lui,
\newblock ``Online robust subspace tracking from partial information,''
\newblock {\em preprint, http://arxiv.org/abs/1109.3827}, 2011.

\end{thebibliography}

\end{document}